\documentclass[12pt, reqno]{amsart}
\usepackage{amsmath, amsthm, amscd, amsfonts, amssymb, graphicx, color}
\usepackage[bookmarksnumbered, colorlinks=false, plainpages]{hyperref}

\newtheorem{theorem}{Theorem}[section]

\newtheorem{proposition}[theorem]{Proposition}

\theoremstyle{definition}
\newtheorem{definition}[theorem]{Definition}
\newtheorem{example}[theorem]{Example}

\theoremstyle{remark}
\newtheorem{remark}[theorem]{Remark}
\numberwithin{equation}{section}

\allowdisplaybreaks
\begin{document}

\title[Controlled $E$-frames in Hilbert spaces]
{Controlled $E$-frames in Hilbert spaces}

\author[H.Hedayatirad]{Hassan Hedayatirad}
\address{Hassan Hedayatirad \\ Department of Mathematics and Computer
Sciences, Hakim Sabzevari University, Sabzevar, P.O. Box 397, IRAN}
\email{ \rm hasan.hedayatirad@hsu.ac.ir; hassan.hedayatirad67@gmail.com}
\author[T.L. Shateri]{Tayebe Lal Shateri }
\address{Tayebe Lal Shateri \\ Department of Mathematics and Computer
Sciences, Hakim Sabzevari University, Sabzevar, P.O. Box 397, IRAN}
\email{ \rm  t.shateri@hsu.ac.ir; t.shateri@gmail.com}
\thanks{*The corresponding author:
t.shateri@hsu.ac.ir; t.shateri@gmail.com (Tayebe Lal Shateri)}
 \subjclass[2010] {Primary 42C15;
Secondary 54D55.} \keywords{ $E$-frame, controlled $E$-frame, Hilbert space, direct sum of Hilbert spaces, dual $E$-frame.}
 \maketitle
\begin{abstract}
In the present paper, we introduce the notion of controlled $E$-frames. Then we investigate and study some properties of them  and characterize all controlled $E$-duals associated with a given controlled $E$-frame. 
\end{abstract}
\section{Introduction}
In 1952, Duffin and Schaeffer \cite{DS} introduced the notion of frames for Hilbert spaces, to study non harmonic Fourier series. In 1986, Daubechies, Grossman, and Meyer \cite{DG} applied the theory of frame to wavelet and Gabor transform. Frames have very important properties which makes them very useful in the characterization of function spaces, signal processing, data compressing, sampling theory and so on. Many authors have done works in this field,  we refer to 
\cite{CAO,CA,CH,HAN} for an introduction to the frame theory and its applications. Various generalizations of frames e.g.
frames of subspaces, wavelet frames, $g$-frames, weighted and controlled frames have developed, see \cite{BAL,CA1,SUN}.

 Controlled frames, as one of the newest generalizations of frames, have been introduced to improve the numerical efficiency of iterative algorithms for inverting the frame operator on abstract Hilbert spaces \cite{BAL}.
Talebi and Dehghan \cite{TAL} introduced the concept of $E$-frames for a separable Hilbert space and study some properties of them.

In the present paper, we introduce the notion of controlled $E$-frames and characterize all controlled $E$-duals associated with a given controlled $E$-frame.

First, we recall some basic notions.

Throughout this paper, we assume that $\mathcal H$ is a separable Hilbert space. A countable family $\lbrace f_{k}\rbrace_{k=1}^{\infty}$ in  $\mathcal{H}$ is called a frame for $\mathcal{H}$ if there exist constants $0<A\leq B<\infty$ such that 
\begin{equation}
A\Vert f\Vert^{2}\leq\sum_{k=1}^{\infty}\big\vert\big<f,f_{k}\big>\big\vert^{2}\leq B\Vert f\Vert^{2}\qquad (f\in\mathcal{H}).
\end{equation}
$A$ and $B$ are called the frame bounds. If just the right inequality holds, we say $\lbrace f_{k}\rbrace_{k=1}^{\infty}$ is a Bessel sequence for $\mathcal{H}$ with bound $B$. A Riesz basis for $\mathcal{H}$ is a family of the form $\lbrace Ue_{k}\rbrace_{k=1}^{\infty}$, where $\lbrace e_{k}\rbrace_{k=1}^{\infty}$ is an orthonormal basis for $\mathcal{H}$ and $U:\mathcal{H}\longrightarrow\mathcal{H}$ is a bounded bijective operator. \\ 
For the sequence $(\mathcal{H}_{n})_{n=1}^{\infty}$ of separable Hilbert spaces, suppose that $$\bigoplus_{n=1}^{\infty}\mathcal{H}_{n}=\bigg\lbrace\lbrace f_{n}\rbrace_{n=1}^{\infty}| f_{n}\in\mathcal{H}_{n},\sum_{n=1}^{\infty}\|f_{n}\|^{2}<\infty\bigg\rbrace.$$ We can define a well defined inner product $\big<.\,,.\big>$ on $\bigoplus_{n=1}^{\infty}\mathcal{H}_{n}$ defined by $$\big<\lbrace f_{n}\rbrace_{n=1}^{\infty},\lbrace g_{n}\rbrace_{n=1}^{\infty}\big>=\sum_{n=1}^{\infty}\big<f_{n},g_{n}\big>.$$ It is well known that $\bigoplus_{n=1}^{\infty}\mathcal{H}_{n}$ is a Hilbert space with respect to this inner product which is called the Hilbert space direct sum of $(\mathcal{H}_{n})_{n=1}^{\infty}$ \cite{Conway}.\\
Let $\mathcal X$ and $\mathcal Y$ be two sequence spaces and $E=(E_{n,k})_{n,k\geq 1}$ be an infinite matrix of real or complex numbers. We say that $E$ defines a matrix mapping from $\mathcal X$ into $\mathcal Y$, if for every sequence $x=\{x_n\}_{n=1}^{\infty}$ in $\mathcal X$, the sequnce $Ex=\{(Ex)_n\}_{n=1}^{\infty}$ is in $\mathcal Y$, where
$$(Ex)_n=\sum_{k=1}^{\infty}E_{n,k}x_k,\quad n=1,2,\ldots.$$
Using the matrix mapping concept,  the authors in \cite{TAL}, have introduced a new notion of frames which is called $E$-frame. Let $E$ be an invertible infinite matrix mapping on $\bigoplus_{n=1}^{\infty}\mathcal{H}$. Then for each $\lbrace f_{k}\rbrace_{k=1}^{\infty}\in\bigoplus_{n=1}^{\infty}\mathcal{H}$, $$E\lbrace f_{k}\rbrace_{k=1}^{\infty}=\bigg\lbrace\sum_{k=1}^{\infty}E_{n,k}f_{k}\bigg\rbrace_{n=1}^{\infty}.$$
\begin{definition}\cite{TAL}
The sequence $\lbrace f_{k}\rbrace_{k=1}^{\infty}$ is called an $E$-frame for $\mathcal{H}$ if there exist constants $0<A\leq B<\infty$ such that 
\begin{equation}\label{$E$-frame}
A\|f\|^{2}\leq\sum_{n=1}^{\infty}\big\vert\big<f,\big(E\lbrace f_{k}\rbrace\big)_{n}\big>\big\vert^{2}\leq B\|f\|^{2}\qquad (f\in\mathcal{H}).
\end{equation}
\end{definition}
If just the right inequality in (\ref{$E$-frame}) holds, then $\lbrace f_{k}\rbrace_{k=1}^{\infty}$ is called an $E$-Bessel sequence for $\mathcal{H}$ with  $E$-Bessel bound $B$. One can easily check that for given constant $B>0$, the sequence $\lbrace f_{k}\rbrace_{k=1}^{\infty}$ is an $E$-Bessel sequence if and only if the operator $T_{E}$ defined by $$T_{E}:\ell^{2}(\mathbb{N}):\longrightarrow\mathcal{H}\,,T\lbrace c_{k}\rbrace_{k=1}^{\infty}=\sum_{n=1}^{\infty}c_{k}\big(E\lbrace f_{k}\rbrace_{k=1}^{\infty}\big)_{n}$$ is a bounded operator from $\ell^{2}(\mathbb{N})$ in to $\mathcal{H}$ with $\|T\|\leq\sqrt{B}$. We call $T_{E}$ the pre $E$-frame operator. Its adjoint, the analysis operator, is given by 
\begin{equation}\label{analysis}
T_{E}^{*}:\mathcal{H}\longrightarrow\ell^{2}(\mathbb{N})\,,T_{E}^{*}f=\bigg\lbrace\big<f,\big(E\lbrace f_{k}\rbrace_{k=1}^{\infty}\big)_{n}\big>\bigg\rbrace_{n=1}^{\infty}.
\end{equation} Composing $T_{E}$ and $T_{E}^{*}$, the $E$-frame operator $$S_{E}:\mathcal{H}\longrightarrow\mathcal{H}\,,S_{E}f=\sum_{n=1}^{\infty}\big<f,\big(E\lbrace f_{k}\rbrace_{k=1}^{\infty}\big)_{n}\big>(E\lbrace f_{k}\rbrace_{k=1}^{\infty}\big)_{n}$$ is obtained. It is easily could be checked that $S_{E}$ is bounded, invertible, self-adjoint and positive \cite{TAL}. This leads us to the following reconstruction formulas
\begin{equation}\label{rec1}
f=\sum_{n=1}^{\infty}\big<f,\big(E\lbrace S_{E}^{-1}f_{k}\rbrace_{k=1}^{\infty}\big)_{n}\big>(E\lbrace f_{k}\rbrace_{k=1}^{\infty}\big)_{n}
\end{equation}
 and $$f=\sum_{n=1}^{\infty}\big<f,\big(E\lbrace f_{k}\rbrace_{k=1}^{\infty}\big)_{n}\big>(E\lbrace S_{E}^{-1}f_{k}\rbrace_{k=1}^{\infty}\big)_{n},$$ for all $f\in\mathcal{H}$.
Let $\lbrace e_{k}\rbrace_{k=1}^{\infty}$ be an orthonormal basis for a separable Hilbert space $\mathcal{H}$. An $E$-Riesz basis for $\mathcal{H}$ is a family of the form $\big\lbrace U\big(E^{-1}\lbrace e_{k}\rbrace_{k=1}^{\infty}\big)_{n}\big\rbrace_{k=1}^{\infty}$, where $U$ is a bounded bijective operator on $\mathcal{H}$.

\section{Controlled $E$-frames}
In this section, we introduce the concept of a controlled $E$-frame in separable Hilbert spaces, and investigate some properties of these frames. We suppose that $\mathcal{H}$ is a separable Hilbert space and $E=(E_{n,k})_{n,k\geq 1}$ is an invertibe infinite matrix mapping on $\bigoplus_{n=1}^{\infty}\mathcal{H}$. We denote by $GL(\mathcal{H})$, the set of all bounded, invertible operators on $\mathcal{H}$ with bounded inverses.
\begin{definition} 
Let $U\in GL(\mathcal{H})$. A sequence $\lbrace \psi_{k}\rbrace_{k=1}^{\infty}$ in $\mathcal{H}$ is called a $U$-controlled $E$-frame if there exist constants $m_{U}>0$ and $M_{U}<\infty$ such that 
\begin{equation}\label{feq}
m_{U}\left\Vert f\right\Vert^{2}\leq\sum_{n=1}^{\infty}\left\langle\left(E\left\lbrace\psi_{k}\right\rbrace_{k=1}^{\infty}\right)_{n},f\right\rangle\left\langle f,U\left(E\left\lbrace\psi_{k}\right\rbrace_{k=1}^{\infty}\right)_{n}\right\rangle\leq M_{U}\left\Vert f\right\Vert^{2}\qquad (f\in H).
\end{equation}
Clearly if $U=id_{\mathcal{H}}$, then the controlled $E$-frame $\lbrace \psi_{k}\rbrace_{k=1}^{\infty}$ is an ordinary $E$-frame for $\mathcal{H}$.
The operator $S_{UE}:\mathcal{H}\longrightarrow\mathcal{H}$ defined by 
\begin{equation}\label{fop}
S_{UE}f=\sum_{n=1}^{\infty}\left\langle\left(E\left\lbrace\psi_{k}\right\rbrace_{k=1}^{\infty}\right)_{n},f\right\rangle U\left(E\left\lbrace\psi_{k}\right\rbrace_{k=1}^{\infty}\right)_{n},
\end{equation}
is called the controlled $E$-frame operator. Using \cite[Proposition 3.2]{BAL} we conclude that $S_{UE}$ is positive, self-adjoint and invertible. Also $S_{UE}=US_{E}$.\\
The operator $T_{UE}:\ell^{2}(\mathbb{N})\longrightarrow\mathcal{H}$ associated with $U$-controlled $E$-frame $\lbrace \psi_{k}\rbrace_{k=1}^{\infty}$ can be defined as 
\begin{equation}
T_{UE}(\lbrace c_{k}\rbrace_{k=1}^{\infty})=\sum_{k=1}^{\infty} c_{k}U\left(E\lbrace \psi_{k}\rbrace_{k=1}^{\infty}\right)_{n}.
\end{equation}
Note that since $U\lbrace\psi_{k}\rbrace_{k=1}^{\infty}$ is a frame for $\mathcal{H}$, so $T_{UE}$ is well defined. We call it the synthesis operator. Now, using (\ref{analysis}) we can represent the controlled $E$-frame operator as $$S_{UE}=T_{UE}T^{*}_{E}.$$
\end{definition}
\begin{theorem}\label{th1}
Let $\lbrace\psi_{k}\rbrace_{k=1}^{\infty}$ be a $U$-controlled $E$-frame for $\mathcal{H}$, then it is an $E$-frame. Moreover, $US_{E}=S_{E}U^{*}$ and so 
\begin{equation*}
\sum_{n=1}^{\infty}\left\langle\left(E\left\lbrace\psi_{k}\right\rbrace_{k=1}^{\infty}\right)_{n},f\right\rangle U\left(E\left\lbrace\psi_{k}\right\rbrace_{k=1}^{\infty}\right)_{n}=\sum_{n=1}^{\infty}\left\langle U\left(E\left\lbrace\psi_{k}\right\rbrace_{k=1}^{\infty}\right)_{n},f\right\rangle\left(E\left\lbrace\psi_{k}\right\rbrace_{k=1}^{\infty}\right)_{n}.
\end{equation*}
\end{theorem}
\begin{proof}
If $\lbrace\psi_{k}\rbrace_{k=1}^{\infty}$ is a $U$-controlled $E$-frame for $\mathcal{H}$, then by (\ref{feq}) we have $$m_{U}\left\Vert f\right\Vert^{2}\leq\left\langle S_{UE}f,f\right\rangle\leq M_{U}\left\Vert f\right\Vert^{2}\qquad (f\in H).$$Hence $m_{U}id_{\mathcal{H}}\leq S_{UE}\leq M_{U}id_{\mathcal{H}}$ and so $S_{UE}\in GL^{+}(\mathcal{H}).$ Now let $\tilde{S}_{E}=U^{-1}S_{UE}$. Then $\tilde{S}\in GL(\mathcal{H})$ and
\begin{align*}
\tilde{S}_{E}f&=U^{-1}S_{UE}f\\&=U^{-1}\left(\sum_{n=1}^{\infty}\left\langle\left(E\left\lbrace\psi_{k}\right\rbrace_{k=1}^{\infty}\right)_{n},f\right\rangle U\left(E\left\lbrace\psi_{k}\right\rbrace_{k=1}^{\infty}\right)_{n}\right)\\&=\sum_{n=1}^{\infty}\left\langle\left(E\left\lbrace\psi_{k}\right\rbrace_{k=1}^{\infty}\right)_{n},f\right\rangle\left(E\left\lbrace\psi_{k}\right\rbrace_{k=1}^{\infty}\right)_{n}=S_{E}f
\end{align*}
for all $f\in\mathcal{H}$. So the $E$-frame operator $S_{E}$ is well defined and $S_{E}=\tilde{S}_{E}\in GL(H)$. Hence $\lbrace \psi_{k}\rbrace_{k=1}^{\infty}$ is an $E$-frame. Clearly $S_{UE}=US_{E}$ and $S_{UE}$ is self adjoint. Thus $S_{UE}^{*}=S_{UE}=US_{E}$ and therefore $US_{E}=S_{UE}=(US_{E})^{*}=S_{E}^{*}U^{*}=S_{E}U^{*}$. 
\end{proof}
Since every controlled $E$-frame is an $E$-frame, (\ref{feq}) yields a criterion to check if a given sequence constitutes an $E$-frame. Furthermore, it becomes obvious from the last result that the role of $U\big(E\lbrace f_{k}\rbrace_{k=1}^{\infty}\big)_{n}$ and $\big(E\lbrace f_{k}\rbrace_{k=1}^{\infty}\big)_{n}$ could have been switched in the definition of controlled frames.
\\\hspace*{10pt} If $U$ is self adjoint we can give necessary and sufficient conditions which are needed for an $E$-frame to form a $U$-controlled $E$-frame: 
\begin{proposition}
Let $U\in GL(\mathcal{H})$ be self adjoint. Then $\lbrace \psi_{k}\rbrace_{k=1}^{\infty}$ is a $U$-controlled $E$-frame for $\mathcal{H}$ if and only if it is an $E$-frame for $\mathcal{H}$ and commutes with the $E$-frame operator $S_{E}$.
\begin{proof}
Suppose that $\lbrace \psi_{k}\rbrace_{k=1}^{\infty}$ is a $U$-controlled $E$-frame. Then by theorem \ref{th1}, $\lbrace \psi_{k}\rbrace_{k=1}^{\infty}$ is an $E$-frame and $US_{E}=S_{E}U^{*}=S_{E}U$. Therefore $S_{UE}S_{E}^{-1}=US_{E}S_{E}^{-1}=S_{E}^{-1}S_{E}U=S_{E}^{-1}US_{E}=S_{E}^{-1}S_{UE}$ and so $U=S_{UE}S_{E}^{-1}$ is positive.\\
Conversely, we note that if $\lbrace \psi_{k}\rbrace_{k=1}^{\infty}$ is an $E$-frame, then $S_{E}\in GL^{+}(\mathcal{H})$. Since $US_{E}=S_{E}U$ and $U$ and $S_{E}$ are positive by assumption, hence $S_{UE}=US_{E}\in GL^{+}(\mathcal{H})$. Now we get the result via \cite[Proposition 2.4]{BAL}
\end{proof} 
\end{proposition}
\begin{theorem}
Let $\lbrace \psi_{k}\rbrace_{k=1}^{\infty}$ be an $E$-frame for $\mathcal{H}$. Then $\lbrace \psi_{k}\rbrace_{k=1}^{\infty}$ is a parseval $U$-controlled $E$-frame for $\mathcal{H}$ if and only if $S_{UE}=US_{E}=id_{\mathcal{H}}$.
\end{theorem}
\begin{proof}
If $\lbrace \psi_{k}\rbrace_{k=1}^{\infty}$ is a parseval controlled $E$-frame, then 
$$\sum_{n=1}^{\infty}\big<f,\big(E\lbrace \psi_{k}\rbrace_{k=1}^{\infty}\big)_{n}\big>\big<U\big(E\lbrace \psi_{k}\rbrace_{k=1}^{\infty}\big)_{n},f\big>=\|f\|^{2}=\big<f,f\big>.$$ Hence $\big<S_{UE}f,f\big>=\big<f,f\big>$ for all $f\in\mathcal{H}$ and so $S_{UE}=US_{E}=id_{\mathcal{H}}$. \\
On the other hand if $S_{UE}=US_{E}=id_{\mathcal{H}}$, then 
\begin{equation*}
\big<f,f\big>=\big<S_{UE}f,f\big>=\sum_{n=1}^{\infty}\big<f,\big(E\lbrace \psi_{k}\rbrace_{k=1}^{\infty}\big)_{n}\big>\big<U\big(E\lbrace \psi_{k}\rbrace_{k=1}^{\infty}\big)_{n},f\big>.
\end{equation*}
\end{proof}
\begin{remark}
If $\lbrace \psi_{k}\rbrace_{k=1}^{\infty}$ is a controlled $E$-frame for $\mathcal{H}$ and $S_{UE}$ is the associated frame operator, then $S_{UE}$ is positive, self adjoint and invertible. So for each $f\in\mathcal{H}$ we have the following reconstruction formula
\begin{align}\label{rec}
f=S_{UE}S_{UE}^{-1}f&=\sum_{n=1}^{\infty}\big<S_{UE}^{-1}f,\big(E\lbrace \psi_{k}\rbrace_{k=1}^{\infty}\big)_{n}\big>U\big(E\lbrace \psi_{k}\rbrace_{k=1}^{\infty}\big)_{n}\\&\nonumber=\sum_{n=1}^{\infty}\big<f,S_{UE}^{-1}\big(E\lbrace \psi_{k}\rbrace_{k=1}^{\infty}\big)_{n}\big>U\big(E\lbrace \psi_{k}\rbrace_{k=1}^{\infty}\big)_{n}
\end{align}
\end{remark}
\begin{definition}
Let $U\in GL(\mathcal{H})$. Suppose that $\lbrace \psi_{k}\rbrace_{k=1}^{\infty}$ is a $U$-controlled $E$-frame and  $\lbrace \phi_{k}\rbrace_{k=1}^{\infty}$ an $E$-Bessel sequence in $\mathcal{H}$. Then we say  $\lbrace \phi_{k}\rbrace_{k=1}^{\infty}$ is a $U$-controlled dual $E$-frame of  $\lbrace \psi_{k}\rbrace_{k=1}^{\infty}$ if
\begin{equation}
f=\sum_{n=1}^{\infty}\big<f,\big(E\lbrace \phi_{k}\rbrace_{k=1}^{\infty}\big)_{n}\big>U\big(E\lbrace \psi_{k}\rbrace_{k=1}^{\infty}\big)_{n}\qquad (f\in\mathcal{H}).
\end{equation}
\end{definition}
In the following example we show that a controlled dual $E$-frame need not be a dual $E$-frame and vice-versa.
\begin{example}
Let $\lbrace e_{k}\rbrace_{k=1}^{\infty}$ be an orthonormal basis for a Hilbert space $\mathcal{H}$ and consider the sequences $\lbrace \psi_{k}\rbrace_{k=1}^{\infty}=\lbrace e_{1},2e_{1},2e_{1}+e_{2},2e_{1}+e_{2}+e_{3},\ldots\rbrace$ and $\lbrace \tilde{\psi}_{k}\rbrace_{k=1}^{\infty}=\lbrace e_{1},2e_{1},2(e_{1}+e_{2}),2(e_{1}+e_{2}+e_{3}),\ldots\rbrace$. Suppose that $E$ is an invertible bi-infinite matrix such that 
\begin{equation*}
E_{n,j}=\begin{cases}
(-1)^{n-j}&n-1\leq j\leq n,\\0&Otherwise.
\end{cases}
\end{equation*}
The matrix form of $E$ is 
\begin{equation*}
E=\begin{pmatrix}
1&0&0&0&0&\cdots\\
-1&1&0&0&0&\cdots\\
0&-1&1&0&0&\cdots\\
0&0&-1&1&0&\cdots\\
0&0&0&-1&1&\cdots\\
\vdots&\vdots&\vdots&\vdots&\ddots&\ddots
\end{pmatrix}.
\end{equation*}
It is easy to show that $\lbrace \psi_{k}\rbrace_{k=1}^{\infty}$ and $\lbrace \tilde{\psi}_{k}\rbrace_{k=1}^{\infty}$ are $E$-frames but they do not constitute a pair of dual $E$-frames. In fact, for given $f\in\mathcal{H}$,
 \begin{align*}
 \sum_{n=1}^{\infty}\left\langle f,\left(E\left\lbrace \psi_{k}\right\rbrace_{k=1}^{\infty}\right)_{n}\right\rangle\left(E\left\lbrace \tilde{\psi}_{k}\right\rbrace_{k=1}^{\infty}\right)_{n}&=\left\langle f,e_{1}\right\rangle e_{1}+\left\langle f,e_{1}\right\rangle e_{1}+2\left\langle f,e_{2}\right\rangle e_{2}+2\left\langle f,e_{3}\right\rangle e_{3}+\cdots\\&=2\sum_{j=1}^{\infty}\left\langle f,e_{j}\right\rangle e_{j}=2f\neq f.
 \end{align*}
 Now consider the operator $U:\mathcal{H}\longrightarrow\mathcal{H}\,;\,U(e_{j})=\dfrac{e_{j}}{2}.$ It is clear that $U\in GL(\mathcal{H})$. We show that $\lbrace \tilde{\psi}_{k}\rbrace_{k=1}^{\infty}$ is a $U$-controlled dual $E$-frame of $\lbrace \psi_{k}\rbrace_{k=1}^{\infty}$. Indeed,
 \begin{align*}
 \sum_{n=1}^{\infty}\left\langle f,\left(E\left\lbrace \psi_{k}\right\rbrace_{k=1}^{\infty}\right)_{n}\right\rangle U\left(E\left\lbrace \tilde{\psi}_{k}\right\rbrace_{k=1}^{\infty}\right)_{n}&=\dfrac{1}{2}\left\langle f,e_{1}\right\rangle e_{1}+\dfrac{1}{2}\left\langle f,e_{1}\right\rangle e_{1}+\dfrac{2}{2}\left\langle f,e_{2}\right\rangle e_{2}+\dfrac{2}{2}\left\langle f,e_{3}\right\rangle e_{3}+\cdots\\&=\sum_{j=1}^{\infty}\left\langle f,e_{j}\right\rangle e_{j}= f.
 \end{align*}
 On the other hand, it is clear that the sequence $\lbrace \phi_{k}\rbrace_{k=1}^{\infty}=\lbrace\dfrac{e_{1}}{3},\dfrac{e_{1}}{3}+\dfrac{2e_{1}}{3},\dfrac{e_{1}}{3}+\dfrac{2e_{1}}{3}+e_{2},\dfrac{e_{1}}{3}+\dfrac{2e_{1}}{3}+e_{2}+e_{3},\cdots\rbrace$ is a dual $E$-frame of $\lbrace \psi_{k}\rbrace_{k=1}^{\infty}$, but it dose not form a $U$-controlled dual $E$-frame for $\lbrace \psi_{k}\rbrace_{k=1}^{\infty}$ because 
\begin{align*}
 \sum_{n=1}^{\infty}\left\langle f,\left(E\left\lbrace \phi_{k}\right\rbrace_{k=1}^{\infty}\right)_{n}\right\rangle U\left(E\left\lbrace\psi_{k}\right\rbrace_{k=1}^{\infty}\right)_{n}&=\dfrac{1}{6}\left\langle f,e_{1}\right\rangle e_{1}+\dfrac{2}{6}\left\langle f,e_{1}\right\rangle e_{1}+\dfrac{1}{2}\left\langle f,e_{2}\right\rangle e_{2}+\dfrac{1}{2}\left\langle f,e_{3}\right\rangle e_{3}+\cdots\\&=\dfrac{1}{2}\sum_{j=1}^{\infty}\left\langle f,e_{j}\right\rangle e_{j}=\dfrac{1}{2}f\neq f.
 \end{align*} 
\end{example}
\begin{definition}\cite{TAL}
Let $\lbrace e_{k}\rbrace_{k=1}^{\infty}$ be an orthonormal basis for $\mathcal{H}$. An $E$-Riesz basis for $\mathcal{H}$ is a family of the form $\left\lbrace V\left(E^{-1}\left\lbrace e_{k}\right\rbrace_{j=1}^{\infty}\right)_{k}\right\rbrace_{k=1}^{\infty}$,where $V$ is a bounded bijective operator on $\mathcal{H}$.
\end{definition} 
\begin{theorem}
Let $\lbrace \psi_{k}\rbrace_{k=1}^{\infty}$ be an $E$-Riesz basis and let $U\in GL(\mathcal{H})$. Then $\lbrace \psi_{k}\rbrace_{k=1}^{\infty}$ is a $U$-controlled $E$-frame for $\mathcal{H}$ if and only if $\lbrace Ve_{j}\rbrace_{j=1}^{\infty}$ is a $U$-controlled frame for $\mathcal{H}$ where $\lbrace e_{j}\rbrace_{j=1}^{\infty}$ is an orthonormal basis for $\mathcal{H}$ and $V$ is a bounded bijection on $\mathcal{H}$ such that $\lbrace \psi_{k}\rbrace_{k=1}^{\infty}=\big\lbrace V\big(E^{-1}\lbrace e_{j}\rbrace_{j=1}^{\infty}\big)_{k}\big\rbrace_{k=1}^{\infty}$.
\end{theorem}
\begin{proof}
Since $\lbrace \psi_{k}\rbrace_{k=1}^{\infty}$ is an $E$-Riesz basis, so there is an orthonormal basis $\lbrace e_{j}\rbrace_{j=1}^{\infty}$ and abounded bijection $V$ on $\mathcal{H}$ such that $\lbrace \psi_{k}\rbrace_{k=1}^{\infty}=\big\lbrace V\big(E^{-1}\lbrace e_{j}\rbrace_{j=1}^{\infty}\big)_{k}\big\rbrace_{k=1}^{\infty}$. Hence for given $f\in\mathcal{H}$ we have
\begin{align*}
&\sum_{n=1}^{\infty}\big<f,\big(E\lbrace \psi_{k}\rbrace_{k=1}^{\infty}\big)_{n}\big>\big<U\big(E\lbrace \psi_{k}\rbrace_{k=1}^{\infty}\big)_{n},f\big>\\&=\sum_{n=1}^{\infty}\big<f,\big(E\big\lbrace V\big(E^{-1}\lbrace e_{j}\rbrace_{j=1}^{\infty}\big)_{k}\big\rbrace_{k=1}^{\infty}\big)_{n}\big>\big<U\big(E\big\lbrace V\big(E^{-1}\lbrace e_{j}\rbrace_{j=1}^{\infty}\big)_{k}\big\rbrace_{k=1}^{\infty}\big)_{n},f\big>\\&=\sum_{n=1}^{\infty}\big<f,V\big(E\big\lbrace\big(E^{-1}\lbrace e_{j}\rbrace_{j=1}^{\infty}\big)_{k}\big\rbrace_{k=1}^{\infty}\big)_{n}\big>\big<UV\big(E\big\lbrace\big(E^{-1}\lbrace e_{j}\rbrace_{j=1}^{\infty}\big)_{k}\big\rbrace_{k=1}^{\infty}\big)_{n},f\big>\\&=\sum_{n=1}^{\infty}\big<f,Ve_{n}\big>\big<UVe_{n},f\big>.
\end{align*}
\end{proof}
\begin{theorem}
Let $U\in GL(\mathcal{H})$ and let $\lbrace \psi_{k}\rbrace_{k=1}^{\infty}$ be a $U$-controlled $E$-frame for $\mathcal{H}$ with the synthesis operator $T_{UE}$. Then $\lbrace \phi_{k}\rbrace_{k=1}^{\infty}$ is a $U$-controlled $E$-dual of $\lbrace \psi_{k}\rbrace_{k=1}^{\infty}$ if and only if $\lbrace \phi_{k}\rbrace_{k=1}^{\infty}=\big\lbrace\big(E^{-1}\lbrace V\delta_{n}\rbrace_{n=1}^{\infty}\big)_{k}\big\rbrace_{k=1}^{\infty}$ where $\lbrace\delta_{n}\rbrace_{n=1}^{\infty}$ is the canonical othonormal basis of $\ell^{2}(\mathbb{N})$ and $V:\ell^{2}(\mathbb{N})\longrightarrow\mathcal{H}$ is a bounded operator such that $T_{UE}V^{*}=id_{\mathcal{H}}$.
\end{theorem}
\begin{proof}
First suppose that $\lbrace \phi_{k}\rbrace_{k=1}^{\infty}$ is a $U$-controlled $E$-dual of $\lbrace \psi_{k}\rbrace_{k=1}^{\infty}$. Then
\begin{equation*}
f=\sum_{n=1}^{\infty}\big<f,\big(E\lbrace \phi_{k}\rbrace_{k=1}^{\infty}\big)_{n}\big>U\big(E\lbrace \psi_{k}\rbrace_{k=1}^{\infty}\big)_{n}\qquad (f\in\mathcal{H}),
\end{equation*}
or equivalently $T_{UE}T_{E}^{*}f=f$ where $T_{E}^{*}$ is the analysis operator of the $E$-Bessel sequence $\lbrace \phi_{k}\rbrace_{k=1}^{\infty}$. Setting $V=T_{E}$ we have $T_{UE}V^{*}=id_{\mathcal{H}}$.\\ 
Conversely, suppose that $V:\ell^{2}(\mathbb{N})\longrightarrow\mathcal{H}$ is a bounded operator such that $T_{UE}V^{*}=id_{\mathcal{H}}$ and let $\phi_k=\big(E^{-1}\lbrace V\delta_{n}\rbrace_{n=1}^{\infty}\big)_{k}$. Then for given $f\in\mathcal{H}$ we have
\begin{align*}
\sum_{k=1}^{\infty}\left\langle f,\left(E\left\lbrace\phi_{k}\right\rbrace_{k=1}^{\infty}\right)_{k}\right\rangle U\left(E\left\lbrace\psi_{k}\right\rbrace_{k=1}^{\infty}\right)_{k}&=\sum_{k=1}^{\infty}\left\langle f,\left(E\left\lbrace\left(E^{-1}\left\lbrace V\delta_{n}\right\rbrace_{n=1}^{\infty}\right)_{j}\right\rbrace_{j=1}^{\infty}\right)_{k}\right\rangle U\left(E\left\lbrace\psi_{k}\right\rbrace_{k=1}^{\infty}\right)_{k}\\&=\sum_{k=1}^{\infty}\left\langle f,V\delta_{k}\right\rangle U\left(E\left\lbrace\psi_{k}\right\rbrace_{k=1}^{\infty}\right)_{k}\\&=\sum_{k=1}^{\infty}\left\langle V^{*}f,\delta_{k}\right\rangle U\left(E\left\lbrace\psi_{k}\right\rbrace_{k=1}^{\infty}\right)_{k}\\&=\sum_{k=1}^{\infty}(V^{*}f)_{k}U\left(E\left\lbrace\psi_{k}\right\rbrace_{k=1}^{\infty}\right)_{k}\\&=T_{UE}V^*f\\&=id_{\mathcal{H}}f=f.
\end{align*}
This proves that $\lbrace \phi_{k}\rbrace_{k=1}^{\infty}$ is a $U$-controlled $E$-dual of $\lbrace \psi_{k}\rbrace_{k=1}^{\infty}$. 
\end{proof}
\begin{theorem}
Let $U\in GL(\mathcal{H})$ and suppose that $\lbrace \psi_{k}\rbrace_{k=1}^{\infty}$ is a $U$-controlled $E$-frame for $\mathcal{H}$ with the synthesis operator $T_{UE}$. Then $\lbrace \phi_{k}\rbrace_{k=1}^{\infty}$ is a $U$-controlled $E$-dual of $\lbrace \psi_{k}\rbrace_{k=1}^{\infty}$ if and only if $\phi_k=S_{UE}^{-1}\psi_k+\big(E^{-1}\lbrace V^{*}\delta_{n}\rbrace_{n=1}^{\infty}\big)_{k}$ where $\lbrace\delta_{n}\rbrace_{n=1}^{\infty}$ is the canonical orthonormal basis of $\ell^{2}(\mathbb{N})$ and $V:\mathcal{H}\longrightarrow\ell^{2}(\mathbb{N})$ is a bounded operator such that $T_{UE}V=0$.
\end{theorem}
\begin{proof}
Suppose that $V:\mathcal{H}\longrightarrow\ell^{2}(\mathbb{N})$ is a bounded operator and $T_{UE}V=0$. Then $\lbrace \phi_{k}\rbrace_{k=1}^{\infty}=\big\lbrace S_{UE}^{-1}\psi_k+\big(E^{-1}\lbrace V^{*}\delta_{n}\rbrace_{n=1}^{\infty}\big)_{k}\big\rbrace_{k=1}^{\infty}$ is an $E$-Bessel sequence. Indeed for given $f\in\mathcal{H}$ we have
\begin{align*}
\sum_{n=1}^{\infty}\left\vert\left\langle f,\left(E\left\lbrace\phi_{k}\right\rbrace_{k=1}^{\infty}\right)_{n}\right\rangle\right\vert^{2}&=\sum_{n=1}^{\infty}\left\vert\left\langle f,\left(E\left\lbrace S_{UE}^{-1}\psi_{k}+\left(E^{-1}\left\lbrace V^{*}\delta_{k}\right\rbrace_{k=1}^{\infty}\right)_{k}\right\rbrace_{k=1}^{\infty}\right)_{n}\right\rangle\right\vert^{2}\\&=\sum_{n=1}^{\infty}\left\vert\left\langle f,\left(E\left\lbrace S_{UE}^{-1}\psi_{k}\right\rbrace_{k=1}^{\infty}\right)_{n}\right\rangle+\left\langle f,\left(E\left\lbrace\left(E^{-1}\left\lbrace V^{*}\delta_{k}\right\rbrace_{k=1}^{\infty}\right)_{k}\right\rbrace_{k=1}^{\infty}\right)_{n}\right\rangle\right\vert^{2}\\&\leq2\left(\sum_{n=1}^{\infty}\left\vert\left\langle f,S_{UE}^{-1}\left(E\left\lbrace\psi_{k}\right\rbrace_{k=1}^{\infty}\right)_{n}\right\rangle\right\vert^{2}+\sum_{n=1}^{\infty}\left\vert\left\langle f,V^{*}\delta_{n}\right\rangle\right\vert^{2}\right)\\&=2\left(\sum_{n=1}^{\infty}\left\vert\left\langle S_{UE}^{-1}f,\left(E\left\lbrace\psi_{k}\right\rbrace_{k=1}^{\infty}\right)_{n}\right\rangle\right\vert^{2}+\left\Vert Vf\right\Vert^{2}\right)\\&\leq 2\left(B\left\Vert S_{UE}^{-1}\right\Vert^{2}+\left\Vert V\right\Vert^{2}\right)\left\Vert f\right\Vert^2,
\end{align*}
where $B$ is an $E$-frame upper bound for $\lbrace \psi_{k}\rbrace_{k=1}^{\infty}$. Moreover
\begin{align*}
\sum_{n=1}^{\infty}\left\langle f,\left(E\left\lbrace\phi_{k}\right\rbrace_{k=1}^{\infty}\right)_{n}\right\rangle U\left(E\left\lbrace\psi_{k}\right\rbrace_{k=1}^{\infty}\right)_{n}&=\sum_{n=1}^{\infty}\left(\left\langle f,\left(E\left\lbrace S_{UE}^{-1}\psi_{k}\right\rbrace_{k=1}^{\infty}\right)_{n}\right\rangle+\left\langle f,V^{*}\delta_{n}\right\rangle\right)U\left(E\left\lbrace\psi_{k}\right\rbrace_{k=1}^{\infty}\right)_{n}\\&=\sum_{n=1}^{\infty}\left\langle f,\left(E\left\lbrace S_{UE}^{-1}\psi_{k}\right\rbrace_{k=1}^{\infty}\right)_{n}\right\rangle U\left(E\left\lbrace\psi_{k}\right\rbrace_{k=1}^{\infty}\right)_{n}\\&+\sum_{n=1}^{\infty}\left\langle f,V^{*}\delta_{n}\right\rangle U\left(E\left\lbrace\psi_{k}\right\rbrace_{k=1}^{\infty}\right)_{n}\\&=f+\sum_{n=1}^{\infty}\left(Vf\right)_{n}U\left(E\left\lbrace\psi_{k}\right\rbrace_{k=1}^{\infty}\right)_{n}\\&=f+T_{UE}Vf=f.
\end{align*}
Therefore $\lbrace \phi_{k}\rbrace_{k=1}^{\infty}$ is a $U$-controlled dual $E$-frame of $\lbrace \psi_{k}\rbrace_{k=1}^{\infty}$.
For the reverse implication let $\lbrace \phi_{k}\rbrace_{k=1}^{\infty}$ be a $U$-controlled dual $E$-frame of $\lbrace \psi_{k}\rbrace_{k=1}^{\infty}$ and define $V=T_\Phi^{*}-T_\Psi^{*}S_{UE}^{-1}$ where $T_\Phi^{*}$ and $T_\Psi^{*}$ are the analysis operators of $E$-frames $\lbrace\phi_{k}\rbrace_{k=1}^{\infty}$ and $\lbrace\psi_{k}\rbrace_{k=1}^{\infty}$, respectively. Clearly $V$ is a bounded operator of $\mathcal{H}$ to $\ell^{2}(\mathbb{N})$. Moreover $$T_{UE}Vf=T_{UE}\big(T_\Phi^{*}-T_\Psi^{*}\big)f=T_{UE}T_\Phi^{*}f-T_{UE}T_\Psi^{*}f=f-f=0.$$
 Finally for each $k\in\mathbb N$, we have
\begin{align*}
S_{UE}^{-1}\phi_k+\left(E^{-1}\left\lbrace V^{*}\delta_{j}\right\rbrace_{j=1}^{\infty}\right)_k&=S_{UE}^{-1}\psi_{k}+\left(E^{-1}\left\lbrace T_{\Phi}\delta_{j}-S_{UE}^{-1}T_{\Psi}\delta_{j}\right\rbrace_{j=1}^{\infty}\right)_{k}\\&=S_{UE}^{-1}\psi_{k}+\left(E^{-1}\left\lbrace\left(E\left\lbrace\phi_{i}\right\rbrace_{i=1}^{\infty}\right)_{j}-S_{UE}^{-1}\left(E\left\lbrace\psi_{i}\right\rbrace_{i=1}^{\infty}\right)_{j}\right\rbrace_{j=1}^{\infty}\right)_{k}\\&=S_{UE}^{-1}\psi_{k}+\left(E^{-1}\left\lbrace\left(E\left\lbrace\phi_{i}\right\rbrace_{i=1}^{\infty}\right)_{j}\right\rbrace_{j=1}^{\infty}\right)_{k}\\&-S_{UE}^{-1}\left(E^{-1}\left\lbrace\left(E\left\lbrace\psi_{i}\right\rbrace_{i=1}^{\infty}\right)_{j}\right\rbrace_{j=1}^{\infty}\right)_{k}\\&=S_{UE}^{-1}\psi_{k}+\phi_{k}-S_{UE}^{-1}\psi_{k}\\&=\phi_{k}.
\end{align*}
This completes the proof.
 \end{proof}
 \begin{theorem}
 Let $\mathcal{H}$ be a Hilbert space and $U\in GL(\mathcal{H})$. Also let $\lbrace \psi_{k}\rbrace_{k=1}^{\infty}$ be a $U$-controlled $E$-frame with synthesis operator $T_{UE}$ and $\lbrace\phi_{k}\rbrace_{k=1}^{\infty}$ be an $E$-Bessel sequence such that 
 \begin{equation}\label{eq5}
 \left\Vert f-\sum_{n=1}^{\infty}\left\langle f,\left(E\left\lbrace\phi_{k}\right\rbrace_{k=1}^{\infty}\right)_{n}\right\rangle U\left(E\left\lbrace\psi_{k}\right\rbrace_{k=1}^{\infty}\right)_{n}\right\Vert<1\qquad(f\in\mathcal{H}).
 \end{equation}
 Then $\left\lbrace\left(D_{E}T_{UE}^{*}\right)^{-1}\phi_{k}\right\rbrace_{k=1}^{\infty}$ is a $U$-controlled dual $E$-frame of $\lbrace\psi_{k}\rbrace_{k=1}^{\infty}$ and 
 \begin{equation*}
\left(D_{E}T_{UE}^{*}\right)^{-1}\phi_{k}=\phi_{k}+\sum_{n=1}^{\infty}\left(id_{\mathcal{H}}-D_{E}T_{UE}^{*}\right)^{n}\phi_{k}
\end{equation*}
where $D_{E}$ is the synthesis $E$-frame operator of $\lbrace\phi_{k}\rbrace_{k=1}^{\infty}$.
 \end{theorem}
 \begin{proof}
 Via \ref{eq5} we can write $\Vert id_{\mathcal{H}}-T_{UE}D_{E}^{*}\Vert<1$, which implies that $T_{UE}D_{E}^{*}$ is an invertible operator. Thus for each $f\in\mathcal{H}$ we have
 \begin{align*}
 f=T_{UE}D_{E}^{*}\left(T_{UE}D_{E}^{*}\right)^{-1}f&=\sum_{n=1}^{\infty}\left\langle\left(T_{UE}D_{E}^{*}\right)^{-1} f,\left(E\left\lbrace\phi_{k}\right\rbrace_{k=1}^{\infty}\right)_{n}\right\rangle U\left(E\left\lbrace\psi_{k}\right\rbrace_{k=1}^{\infty}\right)_{n}\\&=\sum_{n=1}^{\infty}\left\langle f,\left(\left(T_{UE}D_{E}^{*}\right)^{-1}\right)^{*}\left(E\left\lbrace\phi_{k}\right\rbrace_{k=1}^{\infty}\right)_{n}\right\rangle U\left(E\left\lbrace\psi_{k}\right\rbrace_{k=1}^{\infty}\right)_{n}\\&=\sum_{n=1}^{\infty}\left\langle f,\left(D_{E}T_{UE}^{*}\right)^{-1}\left(E\left\lbrace\phi_{k}\right\rbrace_{k=1}^{\infty}\right)_{n}\right\rangle U\left(E\left\lbrace\psi_{k}\right\rbrace_{k=1}^{\infty}\right)_{n}\\&=\sum_{n=1}^{\infty}\left\langle f,\left(E\left\lbrace\left(D_{E}T_{UE}^{*}\right)^{-1}\phi_{k}\right\rbrace_{k=1}^{\infty}\right)_{n}\right\rangle U\left(E\left\lbrace\psi_{k}\right\rbrace_{k=1}^{\infty}\right)_{n}.
 \end{align*}
Hence $\left\lbrace\left(D_{E}T_{UE}^{*}\right)^{-1}\phi_{k}\right\rbrace_{k=1}^{\infty}$ is a $U$-controlled dual $E$-frame of $\lbrace\psi_{k}\rbrace_{k=1}^{\infty}$. Moreover 
we can represent the inverse of $D_{E}T_{UE}^{*}$ by a Neumann series as follow:
\begin{equation*}
\left(D_{E}T_{UE}^{*}\right)^{-1}=\left(id_{\mathcal{H}}-\left(id_{\mathcal{H}}-D_{E}T_{UE}^{*}\right)\right)^{-1}=\sum_{n=0}^{\infty}\left(id_{\mathcal{H}}-D_{E}T_{UE}^{*}\right)^{n}.
\end{equation*}
This completes the proof.
\end{proof}


\end{document}